\documentclass[jfa,final]{elsart}
\usepackage{amssymb}
\newcommand{\Var}{{\rm Var}}

\newcommand{\rf}[1]{(\ref{#1})}
\newcommand{\calV}{{\mathcal V}}
\newcommand{\calA}{{\mathcal A}}
\newcommand{\calB}{{\mathcal B}}
\newcommand{\calF}{{\mathcal F}}

\newcommand{\NC}{{\mathcal N C}}

\newcommand{\TAU}{{\tau\!\!\!\tau}}
\newcommand{\XX}{{\mathbb X}}
\newcommand{\YY}{{\mathbb Y}}
\newcommand{\ZZ}{{\mathbb Z}}
\newcommand{\SSS}{{\mathbb S}}
\newcommand{\UU}{{\mathbb U}}
\newcommand{\VV}{{\mathbb V}}
\newcommand{\WW}{{\mathbb W}}

\newcommand{\II}{{\mathbb I}} %idenitty
\newcommand{\EE}{{\TAU}}

\newcommand{\sR}{{\mathbb R}}
\newcommand{\sC}{{\mathbb C}}

\newcommand{\X}{X}
\newcommand{\Y}{Y}

\newcommand{\E}{{E}}

\newcommand{\mX}{{\mathbf X}}
\newcommand{\mY}{{\mathbf Y}}
\newcommand{\mS}{{\mathbf S}}
\newcommand{\mZ}{{\mathbf Z}}
\newcommand{\mI}{{\mathbf I}}

\newcommand{\la}{\lambda}
\newcommand{\be}{\begin{equation}}
\newcommand{\ee}{\end{equation}}

      \newtheorem{theorem}{Theorem}[section]
       \newtheorem{proposition}[theorem]{Proposition}
       \newtheorem{corollary}[theorem]{Corollary}
       \newtheorem{lemma}[theorem]{Lemma}
  \newtheorem{remark}{Remark}[section]
\begin{document}
\begin{frontmatter}
\title{On a class of free L\'evy laws related to a regression problem}
\runtitle{L\'evy laws and a regression problem}
\newcommand{\email}[1]{\ead{#1}}
\author{Marek Bo\.zejko}
\thanks{\noindent Research partially supported by KBN
grant \#2PO3A00723,  EU Network grant HPRN-CT-2002-00279}

\address{Instytut Matematyczny, Uniwersytet Wroc\l awski, Pl. Grunwaldzki 2/4
50-384 Wroc\l aw, Poland.
}
\email{bozejko@math.uni.wroc.pl}
\author{
W{\l}odzimierz  Bryc
}
\thanks{\noindent Research partially supported by
NSF grants \#INT-0332062, \#DMS-0504198,
and by C.P. Taft Memorial Fund}
%University of Cincinnati's
%Summer Faculty Research Fellowship Program}
%\thanks{\noindent Research partially supported by NSF
%grant \#INT-0332062, by the C.P. Taft Memorial Fund, and University of Cincinnati's
%Summer Faculty Research Fellowship Program}
\address{
Department of Mathematics,
University of Cincinnati,
PO Box 210025,
Cincinnati, OH 45221--0025, USA}
\email{Wlodzimierz.Bryc@UC.edu}

%\date{Printed: \today\  {\tt File: \jobname.TEX}}
%\keywords{}
%\subjclass[2000]{Primary: 46L53; Secondary: 60E05, 05A18}

%\begin{center}
%\fbox{\fbox{\bf\sc Preliminary Draft}}\\ \medskip
%
%\fbox{\bf Not for distribution! Not for citation! Not for reference!}
%\\\medskip \noindent
%\end{center}
%\tolerance=7000

\begin{abstract}
 The
free Meixner laws arise as  the distributions of orthogonal polynomials with
constant-coefficient recursions. We show that these are the laws of the free pairs of random variables
which have linear regressions and quadratic conditional variances when conditioned with respect to their sum.
We apply this result to describe free L\'evy processes with quadratic conditional variances, and to prove a converse
implication related to asymptotic freeness  of random Wishart matrices.
%The $\boxplus$-infinitely divisible  members of this
%family have Wigner's semicircle laws in
%their
%L\'evy-Khinchin representation and appear as the distributions of the free L\'evy processes
%with linear  regressions and quadratic condition variances. The remaining members of the family can be interpreted as
% the "continuous" free analog of the binomial type distributions.
\end{abstract}
%\maketitle
\end{frontmatter}

\section{Introduction}
The family of classical Meixner laws was discovered by Meixner \cite{Meixner-40} who described the class of orthogonal polynomials $p_n(x)$
 with generating function of the form
$\sum_{n=0}^\infty t^n p_n(x)/n!=f(t)e^{x u(t)}$; it turns out that up
to affine transformations of $x$, they correspond to one of the six
classical probability measures: gaussian, Poisson, gamma, Pascal
(negative binomial), binomial, or a two-parameter hyperbolic secant
density which  Schoutens \cite{Schoutens00} calls the (classical)
Meixner law. (Some authors consider only five probability measures
that correspond to non-degenerate polynomials, see
  \cite[page 163]{Chihara}.)
Laha and Lukacs \cite{Laha-Lukacs60} arrived at the same family through the quadratic regression property.
 Morris \cite{Morris82} reinterpreted the latter result into the language of natural exponential families,
and pointed out the connection to orthogonal polynomials.
Weso\l owski \cite{Wesolowski93} proved that the five infinitely-divisible classical Meixner laws
are the laws of  stochastic processes with linear regressions and quadratic conditional variances when the conditional variances
 depend only on
the increments of the process.

The free Meixner systems of polynomials were introduced by
Anshelevich \cite{Anshelevich01} and Saitoh and Yoshida
\cite{Saitoh-Yoshida01}. Combining \cite{Saitoh-Yoshida01} with
\cite[Theorem 4]{Anshelevich01}, up to affine transformations the
corresponding probability measures associated with the free Meixner
system can be classified into six types: Wigner's semicircle, free
Poisson,  free Pascal (free negative binomial), free Gamma,  a law
that we will call pure free Meixner, or a free binomial law. The
first five of these laws are infinitely-divisible with respect to
the additive free convolution and correspond to the  five classical
laws, even though Anshelevich \cite[page 241]{Anshelevich01}
observed that this correspondence does not follow the Bercovici-Pata
bijection \cite{Berkovici-Pata99}. Anshelevich \cite[Remark
6]{Anshelevich01} points out the $q$-Meixner systems of polynomials
that interpolate between the classical case $q=1$ and the free case
$q=0$. The interpolating $q$-deformed Meixner laws appear also as
the laws of classical stochastic processes with linear regressions
and quadratic conditional variances  in \cite{Bryc-Wesolowski-03}.
To facilitate the comparison with the latter paper, in this note we
parameterize the free Meixner family so that our parameters $a,b$
correspond to the parameters $\theta, \tau$ in
\cite{Bryc-Wesolowski-03} when $q=0$
 and the univariate distributions of the processes are taken at time $t=1$.
In Remark \ref{q-interpolation} we observe another
 parallel  between the free and classical cases  on the level of cumulants.

In this note we relate free Meixner laws directly to free probability
via regression properties \rf{LR} and \rf{QV}.
We use the regression characterization
to point out that
the $\boxplus$-infinitely divisible  members of the free Meixner
family  are the distributions of the free L\'evy processes
with linear  regressions and quadratic condition variances, which is a free counterpart of the
 classical result \cite[Theorem 2.1]{Wesolowski89a}. We also prove
  a free analog of the classical characterization \cite[Theorem 2]{Olkin-Rubin-62}
 of  Wishart matrices by the independence of the sum and the quotient.

We conclude this section with the remark that some of the laws in the free Meixner family
 appeared under other names in the literature.
The semicircle law  is the law of the free Brownian motion
\cite[Section 5.3]{Biane98} and appeared as the limiting
distribution of the eigenvalues  of random matrices \cite{Wigner58}.
The free Poisson law \cite[Section 2.7]{Voiculescu00} is also know
as Marchenko-Pastur law \cite{Marchenko-Pastur67}. The free pure
Meixner type  law  appears under the name ``Continuous Binomial" in
\cite[Theorem 4 ]{Anshelevich01}, and for $a=0$ appears as a Brown
measure in \cite[Example 5.6]{Haagerup-Larsen}.  A sub-family of
free Meixner laws appears in random walks on the free group in
Kesten \cite{Kesten}, and as Gaussian and Poisson laws in the
generalized limit theorems \cite{Bozejko-Leinert-Speicher},
\cite{Bozejko-Speicher91}.

%In Section ...
%we give an analytic definition of the free Meixner family, point out the surprising connection with
%the Wigner semicircle laws and with Kesten laws.

\section{Free Meixner laws}
\subsection{Free cumulants}
We assume that our probability space is a von Neumann algebra
$\calA$ with a normal faithful tracial state  $\EE:\calA\to\sC$,
i.e., $\EE(\cdot)$ is linear, continuous in weak* topology,
$\EE(ab)=\EE(ba)$, $\EE(\II)=1$, $\EE(aa^*)\geq 0$ and
$\EE(aa^*)=0$ implies $a=0$. A (noncommutative) random variable
$\XX$ is a self-adjoint ($\XX=\XX^*$) element of $\calA$.

The joint moments of random variables
$\XX_1,\XX_2,\dots,\XX_k\in\calA$
 are complex numbers
$\EE(\XX_1 \XX_2\dots\XX_k)$.
 Since the sequence of univariate
moments $\{\EE(\XX^n):n=0,1,\dots\}$ is real, positive-definite, and
bounded by $\|\XX\|^n$, there is a unique probability measure $\mu$,
called the law of $\XX$, such that $\EE(\XX^n)=\int x^nd\mu$. (In
more algebraic versions of non-commutative probability, a law of
$\XX$ is a unital linear functional $\mu$ on the algebra of
polynomials $\sC(\XX)$ in one variable.)

The concept of freeness was introduced by Voiculescu, see \cite{Voiculescu00} and the references therein.
We are interested in free random variables,
and for our purposes the combinatorial approach of Speicher \cite{Speicher-97}
is convenient.
Recall that a partition $\calV=\{B_1,B_2,\dots\}$ of $\{1,2,\dots,n\}$ is crossing if there are $i_1<j_1<i_2<j_2$ such that
$i_1,i_2$ are in the same block $B_r$ of $\calV$,  $j_1,j_2$ are in the same block $B_s\in\calV$ and $B_r\ne B_s$.
By $\NC(n)$ we denote the set of all non-crossing partitions of $\{1,2,\dots,n\}$.

Let $\sC\langle\XX_1,\dots,\XX_k\rangle$ denote the non-commutative
polynomials in variables $\XX_1,\dots,\XX_k$.
The free (non-crossing) cumulants %$R_k$
are the $k$-linear maps $R_k:\sC\langle\XX_1,\dots,\XX_k\rangle\to\sC$
defined recurrently by the formula that connects them with moments:
\begin{equation}\label{Cumulants}
\EE(\XX_1 \XX_2\dots\XX_n)=\sum_{\calV\in\NC(n)} R_\calV(\XX_1, \XX_2,\dots,\XX_n),
\end{equation}
where
\[
R_\calV(\XX_1, \XX_2,\dots,\XX_n)=\prod_{B\in\calV}R_{|B|}(\XX_j:j\in B),
\]
see \cite{Speicher-97}. We will also write $R_n(\XX)$ for $R_n(\XX,\XX,\dots,\XX)$, and $R_n(\mu)$ for $R_n(\XX)$ when
$\mu$ is the law of $\XX$.

Random variables $\XX_1,\XX_2,\dots,\XX_k\in\calA$ are free if for every $n\geq 1$ and every
non-constant choice of
$\ZZ_1,\dots,\ZZ_n\in\{\XX_1,\XX_2,\dots,\XX_k\}$
we have
\[R_n(\ZZ_1,\ZZ_2,\dots,\ZZ_n)=0,
\]
see \cite{Speicher-97}. The free convolution $\mu\boxplus\nu$ of measures $\mu,\nu$ is the law
of $\XX+\YY$ where $\XX,\YY$ are free and have laws $\mu,\nu$ respectively.

The $R$-transform of a random variable $\XX$ is
\[r_\XX(z)=\sum_{k=0}^\infty R_{k+1}(\XX)z^k.\]
It  linearizes the additive free convolution, $r_{\mu\boxplus\nu}(z)=r_\mu(z)+r_\nu(z)$, see
\cite{VDN92}, \cite{Voiculescu00}, or \cite{Hiai-Petz}.

%\subsection{Free infinitely divisible probability measures}
%$$R_aX(z)=a R_x(az)$$

\subsection{Free Meixner laws}
We are interested in the two-parameter family of probability
measures $\{\mu_{a,b}: a\in\sR, b\geq -1\}$ with the
Cauchy-Stieltjes transform
\begin{equation}
\label{free-Meixner}
\int_{\sR}\frac{1}{z-y}\mu_{a,b}(dy)=
\frac{ (1+2b)z+a  -
    \sqrt{ (z-a)^2 -
        4 (1+b) }}{2(b z^2+a z+1 )},
\end{equation}
which we will call the free (standardized, i.e. with mean zero and  variance one) Meixner laws.
The absolutely continuous part of  $\mu_{a,b}$ is
\[\frac{\sqrt{4(1+b)-(x-a)^2}}{2\pi (bx^2 +ax +1)}\]
 on $a-2\sqrt{1+b}\leq x\leq a+2\sqrt{1+b}$; the measure may also have one atom if $a^2>4b\geq 0$,
and a second atom if $-1\leq b<0$. For more details, see  \cite{Saitoh-Yoshida01}, who
analyze
  monic orthogonal
polynomials which satisfy constant-coefficient recurrences; in our parametrization
the monic orthogonal polynomials  $p_n$ with respect to $\mu_{a,b}$ satisfy the recurrence
\begin{equation}\label{ccc}
xp_n(x)=p_{n+1}(x)+a p_n(x)+(1+b) p_{n-1}(x), \, n\geq2
\end{equation}
with the perturbed initial terms $p_0(x)=1$, $p_1(x)=x$, $p_2(x)=x^2-ax-1$.
This can be seen
 from the corresponding continued fraction representation
\[\displaystyle
\int_{\sR}\frac{1}{z-y}\mu_{a,b}(dy)=\displaystyle
\frac{1}{\displaystyle z-\frac{\displaystyle  1}{z-a-
\frac{\displaystyle 1+b}{\displaystyle z-a-\frac{1+b}{\ddots}}}}
.\]

For ease of reference we state the corresponding $R$-transform
\begin{equation}\label{R0}
r_{\mu_{a,b}}(z)=\frac{2z}{1-az+\sqrt{(1-az)^2-4z^2b}}.
\end{equation}

Since \rf{free-Meixner} defines measures with mean zero and variance 1,
 it is convenient to enlarge this class by allowing
translations by $t$ and dilations by $r$; that is, to consider probability laws up to their type.
Recall that a $\mu$-type law is a law of an arbitrary affine
 function of a random variable with the law $\mu$. In other words, the
  $\mu$-type laws consist of all probability laws $\{D_r(\mu)\boxplus \delta_t: r, t\in\sR\}$, where
  the dilation $D_r$ is defined as $D_r(\mu)(A)=\mu(A/r)$ if $r\ne0$ and $D_0(\mu)=\delta_0$.

 The laws of free Meixner type correspond to various
 reparametrizations of \rf{free-Meixner}
 and occurred in many places in the literature, see \cite{Anshelevich01}, \cite{Bozejko-Leinert-Speicher},
  \cite{Bozejko-Speicher91},  \cite[Corollary 7.2]{Capitaine-Casalis04},
 \cite[Example 5.6]{Haagerup-Larsen}, \cite{Kesten}, \cite{McKay}, \cite{Saitoh-Yoshida01}.
 In particular, recursion \rf{ccc} for $b\geq 0$ is
 a reparametrization of five recursions that
 appear in Anshelevich \cite[Theorem 4]{Anshelevich01}.

%\input{meixner.skt}

%The following is \cite[Theorem 3.2(ii)]{Saitoh-Yoshida01} in our parametrization.
%\begin{proposition}[Saitoh-Yoshida]\label{P1} Let $\mu_{a,b}$ be the (standardized) free Meixner distribution with $b\geq0$.
%Then $\mu_{a,b}$ is $\boxplus$-infinitely divisible and  its
% L\'evy-Khinchin measure
%$\nu$
%\rf{Levy-Khinchin} is the Wigner semicircle law $\omega_{a,b}$ of mean $a$ and variance $b$.
%\end{proposition}

%\begin{pf}
%The Cauchy-Stieltjes transform $G_{\omega_{a,b}}(z)=\int_{\sR}1/(z-y){\omega_{a,b}}(dx)$ of the semicircle distribution is
%$$
%G_{\omega_{a,b}}(z)=\frac{2}{z-a+\sqrt{(a-z)^2-4b}}
%%\frac{z-a - {\sqrt{ {\left( a - z \right) }^2-4 b }}}{2 b}
%$$
%and the corresponding moment generating function is $$M_{\omega_{a,b}}(z)=G_{\omega_{a,b}}(1/z)/z=r_{\mu_{a,b}}(z)/z.$$
%Thus $R_1(\mu_{a,b})=0$ and the higher cumulants  are
% $R_{n+2}(\mu_{a,b})=\int x^{n}d{\omega_{a,b}}(x)$, identifying the L\'evy-Khinchin measure \rf{Levy-Khinchin} through \rf{RRR}.
%\end{pf}

Saitoh and Yoshida \cite[Theorem 3.2]{Saitoh-Yoshida01} transcribed into our notation says
that for $b\geq 0$ (standardized) free Meixner law  $\mu_{a,b}$ is $\boxplus$-infinitely divisible, i.e.,
for every integer $n$ there exists a measure $\nu_n$ such that
$\mu_{a,b}=\nu_n^{\boxplus n}$, and that the corresponding
 L\'evy-Khinchin measure
\rf{Levy-Khinchin} is the Wigner semicircle law $\omega_{a,b}$ of mean $a$ and variance $b$, i.e.
the $R$ transform of $\mu_{a,b}$ is given by
\begin{equation}\label{Levy-Khinchin}
r_{\mu_{a,b}}(z)=\int \frac{z}{1-xz}\omega_{a,b}(dx);
\end{equation}
see also \cite[page 242]{Anshelevich01}. (Here we follow the  L\'evy-Khinchin representation
from \cite[Theorem 3.3.6]{Hiai-Petz};
other authors state the L\'evy-Khinchin formula in another form, see
\cite{Berkovici-Voiculescu92}, or \cite{Barndorff-Nielsen-Thorbj02c}.
The "L\'evy-Khinchin"  measures that enter these two representations differ by a factor of $(1+x^2)$, see  \cite[Theorem 2.7]{Barndorff-Nielsen-Thorbj02a} or
\cite[Section 6.1]{Anshelevich01a} who also discusses
the $q$-interpolation of  \rf{Levy-Khinchin}.)
From \rf{Levy-Khinchin} it follows that  the free cumulants of  $\mu_{a,b}$ are $R_1(\mu_{a,b})=0$ and
\begin{equation}\label{RRR}
R_{n+2}(\mu_{a,b})=\int x^n \omega_{a,b}(dx), \, n\geq 0.
\end{equation}
%The latter is the basis for $q$-interpolation in \cite[Section 6]{Anshelevich01a}.

As observed in \cite[Theorem 3.2(i)]{Saitoh-Yoshida01}, the free Meixner type laws with $-1\leq b<0$ are not
$\boxplus$-infinitely divisible. In fact, they are the free counterpart of the classical binomial laws
$\left((1-p)\delta_0+p\delta_1\right)^{*n}$. This was first noticed for symmetric free Meixner laws in
\cite{Bozejko-Wysoczanski01}, where the authors prove that
\[\mu_{0,b}=D_{1/\sqrt{t}}\left(\frac12\delta_{-1}+\frac{1}{2}\delta_1\right)^{\boxplus t}\]
with $t=-1/b\geq 1$. Compare also \cite[Example 3.4]{Saitoh-Yoshida01}.
(The
additive $t$-fold free convolution $\mu^{\boxplus t}$
is well defined for the continuous range of values  $t\geq 1$, see
 \cite{Nica-Speicher}.)
Since
\[\mu_{a,-1}
=\frac{1}{2}\left(1+
\frac{a}{\sqrt{1+a^2}}\right)\delta_{\frac{1}{2}(a-\sqrt{4+a^2})}+
\frac{1}{2}\left(1-\frac{a}{\sqrt{1+a^2}}\right)\delta_{\frac{1}{2}(a+\sqrt{4+a^2})}
\]
is supported on two-points, the following is a generalization  of the above result.
\begin{proposition}\label{P2} If $\mu_{a,b}$ is a free Meixner measure \rf{free-Meixner} with parameters $a\in\sR$, $-1\leq b<0$ and $t=-1/b$ then
\begin{equation}\label{kesten}
\mu_{a,b}=D_{\sqrt{|b|}}
\left(
\mu_{a/\sqrt{|b|},-1}
^{\boxplus t}\right).
\end{equation}
\end{proposition}
\begin{pf}
We observe in more generality that if $\mu_{a,b}$ is given by \rf{free-Meixner} then
\[
D_{1/\la}(\mu_{a,b}^{\boxplus \la^2})=\mu_{a/\la, b/\la^2}, \; \la\ne 0.
\]
This follows from  \rf{R0} and the fact that the $R$ transform of the dilatation
$D_\la(\mu)$ is $\la r_{\mu}(\la z)$.
%
%
%Put $t=-1/b $ and
%$$ p_1=\frac{1+ a/\sqrt{a^2+4|b|}}{2}, \;p_2=1-p_1,$$
%$$x_1=\frac{a-\sqrt{a^2+4|b|}}{2\sqrt{|b|}},$$
%$$x_2=\frac{a+\sqrt{a^2+4|b|}}{2\sqrt{|b|}}.$$
%%
%%\fbox{Dylatacje?}
%%$$\nu_{a,b}=\frac{1}{2}(1+ a/\sqrt{1+a^2})\delta_{\frac{1}{2}(a-\sqrt{4+a^2})}+
%%\frac{1}{2}(1-a/\sqrt{1+a^2})\delta_{\frac{1}{2}(a+\sqrt{4+a^2})}$$
%%%$$\nu_{a,b}=\frac{1}{2}(1+ a/\sqrt{b+a^2})\delta_{x_1}+
%%%\frac{1}{2}(1-a/\sqrt{b+a^2})\delta_{x_2}$$
%%denote the two-point discrete measure.
%Then $\mu_{a,-1}=p_1\delta_{x_1}+p_2\delta_{x_2}$. A calculation gives the $R$ transform of $\mu_{a,-1}$ as
%$$
%r_{\mu_{a,-1}}(z)=\frac{z-a/\sqrt{|b|}}{z^2-1-za/\sqrt{|b|}}.
%$$
%Therefore, the $R$ transform of $\mu_{a,-1}^{\boxplus t}$ is
%$$
%r_{\mu_{a,-1}^{\boxplus t}}(z)=\frac{z-a/\sqrt{|b|}}{(1+za/\sqrt{|b|}-z^2)b}.
%$$
%The dilatation by $\sqrt{|b|}$ acts on the $R$ transform via
%$$r_{D(\nu^{\boxplus t})}(z)=\sqrt{|b|}r_{\nu^{\boxplus t}}(z\sqrt{|b|}),$$
%we get
%
%\fbox{to be continued...}
\end{pf}
%We remark that the proof shows that  all free Meixner type measures  are generated by dilations, translations, and
%additive free convolution powers from a one-parameter family $\{\mu_{a,\pm 1}:a\in\sR\}$.

The similarity of $\mu_{a,b}$ with $b<0$ to the classical binomial law is further strengthened by
comparing Theorem \ref{Laha-Lukacs}(vi) with  Theorem \ref{T1}(vi), and by random projection asymptotic of
\cite[Section 3.2]{Collins04}.

Next we show that the free cumulants of the free Meixner law  can be expressed as sums over the
 non-crossing partitions $\NC_{\leq 2}(n)$ with blocks of size at most $2$.
 \begin{proposition} If $\mu_{a,b}$ is a (standardized) free Meixner law with parameters $a\in\sR,b\geq -1$,
 then the free cumulants are $R_1(\mu_{a,b})=0$ and
\begin{equation}\label{CC}
 R_{n+2}(\mu_{a,b})=\sum_{\calV\in\NC_{\leq 2}(n)}a^{s(\calV)}b^{|\calV|-s(\calV)}, \; n\geq 0,
\end{equation}
where  $s(\calV)$ is the number of blocks of size $1$ (singletons)
of a partition $\calV\in\NC_{\leq 2}(n)$.
\end{proposition}
For example
$R_5(\mu_{a,b})=a^3+3ab$,  see Figure \ref{Fig 2}.
% Similarly $R_6=a^4+6 a^2b+ 2 b^2$,
%$R_7=a^5+10 a^3 b+10 a b^2$,
%$R_8=a^6+15a^4 b+30 a^2b^2+5 b^3$, $R_9=a^7+21 a^5 b+70 a^3b^3+35 ab^3$.
\begin{figure}[hbt]
\begin{picture}(100,100)(-120,0)

\put(0, 40){\circle*{4}}
\put(10, 40){\circle*{4}}
\put(20, 40){\circle*{4}}
\put(10,40){\makebox(0,-10)[ct]{$\calV_1$}}
\put(10,80){\makebox(0,5)[ct]{$a^3$}}

\put(40, 40){\circle*{4}}
\put(50, 40){\circle*{4}}
\put(60, 40){\circle*{4}}

\put(50,40){\makebox(0,-10)[ct]{$\calV_2$}}
\put(50,80){\makebox(0,5)[ct]{$ab$}}
\put(100, 40){\circle*{4}}
\put(90, 40){\circle*{4}}
\put(80, 40){\circle*{4}}
\put(90,40){\makebox(0,-10)[ct]{$\calV_3$}}
\put(90,80){\makebox(0,5)[ct]{$ab$}}

\put(120, 40){\circle*{4}}
\put(130, 40){\circle*{4}}
\put(140, 40){\circle*{4}}
\put(130,40){\makebox(0,-10)[ct]{$\calV_4$}}
\put(130,80){\makebox(0,5)[ct]{$ab$}}

\thicklines
\put(0,40){\line(0,1){20}}
\put(10,40){\line(0,1){20}}
\put(20,40){\line(0,1){20}}

\put(40,40){\line(0,1){20}}
\put(50,40){\line(0,1){20}}
\put(60,40){\line(0,1){20}}

\put(80,40){\line(0,1){20}}
\put(90,40){\line(0,1){20}}
\put(100,40){\line(0,1){20}}

\put(120,40){\line(0,1){20}}
\put(130,40){\line(0,1){15}}
\put(140,40){\line(0,1){20}}

\put(120,60){\line(1,0){20}}

\put(80,60){\line(1,0){10}}

\put(50,60){\line(1,0){10}}
%\put(60,60){\line(1,0){5}}
%\put(20,42){\line(0,1){10}}\put(30,42){\line(0,1){8}}\put(50,42){\line(0,1){10}}
%\put(20,52){\line(1,0){5}}%\put(30,50){\line(1,0){5}}
%\put(50,52){\line(-1,0){5}}
%\put(10,60){\line(1,0){60}}
\end{picture}
\caption{$\NC_{\leq 2}(3)=\left\{\calV_1,\calV_2,\calV_3,\calV_4\right\}$. \label{Fig 2}}
\end{figure}

\begin{pf}
For $b\geq 0$ the monic orthogonal polynomials $q_n(x)$ corresponding to
the semicircle law $\omega_{a,b}$ of mean $a$ and variance $b$
 satisfy the three step recursion
\[
(x-a)q_n(x)=q_{n+1}(x)+b q_{n-1}(x).
\]
Since by \cite[Theorem 3.2]{Saitoh-Yoshida01}, semicircle law
 $\omega_{a,b}$ is the L\'evy-Khinchin measure
 \rf{RRR} for $\mu_{a,b}$, from  \cite[Corollary 5.1]{Accardi-Bozejko98} we deduce \rf{CC}.
 Since the
 cumulants
 are given by  algebraic expressions,  formula \rf{CC} extends to all $-1\leq b<\infty$.
\end{pf}

\section{Conditional moments and free Meixner laws}
We begin with the classical characterization which, as observed in \cite{Wesolowski89a}, follows from the argument in
\cite{Laha-Lukacs60}.
\begin{theorem} \label{Laha-Lukacs}
Suppose $X,Y$ are non-degenerate
independent square-integrable classical random variables,
and there are numbers $\alpha,\alpha_0,C,a,b \in\sR$ such that
  \begin{equation}\label{LR0}
%\[
\E(\X|\X+\Y)=\alpha (\X+\Y)+\alpha_0
%\]
\end{equation}
 and %there are numbers $C,a,b,  \in\sR$ such that
\begin{equation}\label{QV0}
%\[
\Var(\X|\X+\Y)=C(1+a(\X+\Y)+b(\X+\Y)^2).
%\]
\end{equation}
Then $\X$ and $\Y$   have the classical Meixner type laws.
In particular, if  $\E(\X)=\E(\Y)=0$,
$\E(\X^2)+\E(\Y^2)=1$ then $\alpha=\E(\X^2)$, $\alpha_0=0$, $C=\alpha(1-\alpha)/(1+b)$,
and the law of $\X$ is
\begin{enumerate}
\item normal (gaussian), if $a=b=0$;
\item Poisson type, if $b=0$ and $a\ne 0$;
\item Pascal type, if $b>0$ and $a^2>4b$;
\item Gamma type, if $b>0$ and $a^2=4b$;
\item Meixner type, if $b>0$ and $a^2<4b$;
\item Binomial type, if $b=-\alpha/n=-(1-\alpha)/m$ and $m,n$ are integers.
\end{enumerate}
\end{theorem}
%\begin{pf} For $|t|$ small enough, let
%$r_X(t)$, $r_Y(t)$ be the derivatives of the
%logarithms of the characteristic functions of $X,Y$ respectively, compare \cite{Nica-96}. It is well known, see
%\cite{Kagan-Linnik-Rao} that the assumptions on conditional moments imply
%$$
%(1-c)r_Y(t)=cr_X(t)
%$$
%and
%$$
%r_X'(t)=
%$$
%\end{pf}
%\renewcommand{\thetheorem}{\arabic{section}.\arabic{theorem}}

To state the free version of this theorem, recall that
if $\calB\subset\calA$ is a von Neumann subalgebra of a von Neumann algebra $\calA$ with a normalized trace $\EE$, then
 there exists a unique conditional expectation
from $\calA$ to $\calB$ with respect to $\EE$, see  \cite[Vol. I page 332]{TakesakiI-III}, which we denote by
$\EE(\cdot|\calB)$; the conditional expectation of a self-adjoint element $\XX\in\calA$
is a unique self-adjoint element of $\calB$. %% Added for reviewer

The conditional variance is as usual
$\Var(\XX|\calB)=\EE(\XX^2|\calB)-\left(\EE(\XX|\calB)\right)^2$. If $\YY=\EE(\XX|\calB)$ then $\YY\in\calB$ so
$\EE(\XX\YY|\calB)=\EE(\XX|\calB)\YY=\YY^2$ and similarly $\EE(\YY\XX|\calB)=\YY^2$. Thus
\begin{equation}\label{Var2}
\Var(\XX|\calB)=\EE\left.\left((\XX-\YY)^2\right|\calB\right).
\end{equation}

For fixed $\XX\in\calA$ by $\EE(\cdot|{\XX})$ we denote the conditional expectation corresponding to the von Neumann
algebra $\calB$ generated by $\XX$. A random variable $\XX$ is non-degenerate if it is not a multiple of identity; under faithful
state this is equivalent to $\Var(\XX)\ne0$.

\begin{theorem} \label{T1} Suppose $\XX,\YY\in\calA$ are free, self-adjoint, non-degenerate
and there are numbers $\alpha,\alpha_0,C,a,b \in\sR$ such that
  \begin{equation}\label{LR}
\EE(\XX|\XX+\YY)=\alpha (\XX+\YY)+\alpha_0\II
\end{equation}
 and %there are numbers $C,a,b,  \in\sR$ such that
\begin{equation}\label{QV}
\Var(\XX|\XX+\YY)=C(\II+a(\XX+\YY)+b(\XX+\YY)^2).
\end{equation}
Then $\XX$ and $\YY$   have the free Meixner type laws.
In particular, if  $\EE(\XX)=\EE(\YY)=0$,
$\EE(\XX^2)+\EE(\YY^2)=1$ then $\alpha=\EE(\XX^2)$, $\alpha_0=0$, $C=\alpha(1-\alpha)/(1+b)$, $b\geq -\min\{\alpha,1-\alpha\}$,
and
$\XX/\sqrt{\alpha}$ has the
 $\mu_{a/\sqrt{\alpha},b/\alpha}$ law with the Cauchy-Stieltjes transform \rf{free-Meixner}.
In particular, the law of $\XX$ is
\begin{enumerate}
\item the Wigner's semicircle law if $a=b=0$;
\item  the free Poisson type law if $b=0$ and $a\ne 0$;
\item the free Pascal (negative binomial) type law if $b>0$ and $a^2>4b$;
\item the free Gamma type law if $b>0$ and $a^2=4b$;
\item the pure free Meixner type law if $b> 0$ and $a^2<4b$;
\item the free binomial type law  \rf{kesten} if $-\min\{\alpha,1-\alpha\}\leq b <0$.
\end{enumerate}
\end{theorem}

\input{parabola.skt}

The proof of this theorem is given in Section \ref{proof of T1}.
We now list some consequences.

We need the following properties of conditional expectations.
\begin{lemma}\label{Conditionings}\begin{enumerate}
\item If $\XX\in\calA,\YY\in\calB$, then
\begin{equation}\label{conditioning}
\EE(\XX\YY)=\EE(\EE(\XX|\calB)\YY)
\end{equation}
\item If random variables $\UU,\VV\in\calA$ are free then $\EE(\UU|\VV)=\EE(\UU)\II$.
\item Let $\WW$ be a (self-adjoint) element of the von Neumann algebra generated by a self-adjoint $\VV\in\calA$.
If for all $n\geq 1$ we have
$\EE(\UU\VV^n)=\EE(\WW\VV^n)$
then
$
\EE(\UU|\VV)=\WW$.
\item If $\EE(\UU_1\VV^n)=\EE(\UU_2\VV^n)$ for all $n\geq 1$,
then
$\EE(\UU_1|\VV)=\EE(\UU_2|\VV)$.
\item If $\UU$ is self-adjoint then $\EE((\EE(\UU|\calB))^2)\leq \EE(\UU^2)$.
\end{enumerate}
\end{lemma}
\begin{pf} (i) This follows from the definition, see \cite{Takesaki} or \cite[Vol. II page 211]{TakesakiI-III}.
(ii) If $\ZZ$ is in the von Neumann algebra generated by $\VV$, then
$\EE((\UU-c\II)\ZZ)=\EE(\UU-c\II)\EE(\ZZ)$. Applying this to
$\ZZ=\EE(\UU|\VV)-\EE(\UU)\II$ and $c=\EE(\UU)$ after taking into account \rf{conditioning} we get
$
\EE(\ZZ^2)=\EE\left(\ZZ (\EE(\UU|\VV)-c\II)\right)=\EE\left(\ZZ \EE(\UU-c\II|\VV)\right)=
\EE(\ZZ (\UU-c\II))=\EE(\ZZ)\EE(\UU-c\II)=0
$.
Thus $\EE(\UU|\VV)=\EE(\UU)\II$.

(iii) Let $\WW'=\EE(\UU|\VV)$. Then $\EE\left((\WW-\WW')p(\VV)\right)=0$ for all polynomials $p$.
Since polynomials $p(\VV)$ are  dense in the von Neumann algebra generated by $\VV$, and $\EE(\cdot)$ is normal,
this implies that
$\EE((\WW-\WW')(\WW-\WW')^*)=0$; by faithfulness of $\EE(\cdot)$ we deduce that
$\WW'=\WW$.

(iv) Apply (iii). (v) See \cite[Vol II page 211]{TakesakiI-III}.

% Since $\EE(\UU|\calB)\in\calB$, we have  $\EE(\UU|\calB)^2=\EE(\EE(\UU|\calB)\UU|\calB)$, so the expectations
%are real and by (non-commutative) Cauchy-Schwartz inequality
%$$
%\EE\left(\EE(\UU|\calB)^2\right)=\EE\left(\UU\EE(\UU|\calB)\right)\leq \left(\EE\left(\EE(\UU|\calB)^2\right)\right)^{1/2}\left(\EE(\UU^2)\right)^{1/2}
%$$
\end{pf}

Recall that a non-commutative stochastic
process $(\XX_t)_{t\geq 0}$ on a probability space $(\calA,\EE)$ is a mapping $[0,\infty)\to\calA$.
A stochastic process $(\XX_t)$ has (additive) free increments if
for every $t_1<t_2<\dots<t_k$ random variables
\[
\XX_{t_1},\XX_{t_2}-\XX_{t_1},\dots,\XX_{t_{k}}-\XX_{t_{k-1}}
\]
are free. A stochastic process $(\XX_t)$ is a free L\'evy process if the law of $\XX_t$ converges to $\delta_0$ as $t\to 0$ and
it has free and stationary increments, i.e.
for $s<t$,  the law of
$\XX_t-\XX_s$ is the same as the law of $\XX_{t-s}$.   (Compare \cite[Section 5]{Barndorff-Nielsen-Thorbj02a}.)

We remark that martingale properties of  free L\'evy processes follow from the classical arguments.
In particular, if $(\XX_t)$ is a free L\'evy process such that
$\EE(\XX_t)=0$ and $\EE(\XX_t^2)=t$ for all $t>0$ then
\begin{equation}\label{Rev Mart}
\EE(\XX_s|\XX_u)=\frac{s}{u} \XX_u
\end{equation}
for all $s<u$.
For completeness we include the proof. Suppose $s=\frac{m}{n} u$ for some integers $m<n$. Let
$$\SSS_k=\sum_{j=1}^k(\XX_{j u/n}-\XX_{(j-1)u/n}).$$
 Then $\SSS_k$ is the sum of free random variables with the same distribution,
  and by exchangeability (Lemma \ref{Conditionings}(iii)) we have
$$
\EE(\XX_s|\XX_u)=\EE(\SSS_m|\SSS_n)=m\EE(\SSS_1|\SSS_n)=\frac{m}{n}\SSS_n=\frac{s}{u}\XX_u.
$$
Suppose now $s<u$ is arbitrary. Take a sequence $s_q\to s$ of the previous form. Since
$\EE(\XX_s-\XX_{s_q})^2=|s-s_q|\to 0$, by Lemma \ref{Conditionings}(v) we get \rf{Rev Mart}.

The following is a free version of two different classical probability results:
 \cite[Theorem 2.1]{Wesolowski89a} which characterizes classical L\'evy processes with quadratic conditional variances,
 and  \cite[Theorem 4.3]{Bryc-Wesolowski-03} which characterizes intrinsically the
  classical versions of the free L\'evy processes
that satisfy \rf{QV'}.
\begin{proposition}\label{SP} Suppose $(\XX_t)_{t\geq 0}$ is a free L\'evy process such that
 $\EE(\XX_t)=0,\EE(\XX_t^2)=t$, and that there are constants $\alpha,\beta\in\sR$ and a normalizing function $C_{t,u}$ such that
for all $t<u$
%$$
%\EE(\XX_t|\XX_u)=\frac{t}{u}\XX_u
%$$
\begin{equation}\label{QV'}
%\EE\left.\left(\XX_t^2\right|\XX_u\right)=\frac{t^2}{u^2}\XX_u^2+C(t,u)\left(1+\frac{\alpha}{u}\XX_u+\frac{\beta}{u^2}\XX_u^2\right).
\Var(\XX_t|\XX_u)=C_{t,u}\left(1+\frac{\eta}{u}\XX_u+\frac{\sigma}{u^2}\XX_u^2\right).
\end{equation}
Then for $s,t>0$ the increment   $\XX_{s+t}-\XX_s$ has the
free Meixner type  law  $\mu_{a,b}$ with parameters $a=\eta/\sqrt{t}\in\sR, b=\sigma/t\geq0$.
Moreover, the $R$ transform of $\XX_t$ is
%$$
%r_{X_t}(z)=t \frac{1 - \alpha z -\sqrt{(1-\alpha z)^2-4 z^2\beta}}{2\beta z}
%$$
\[
r_{\XX_t}(z)=\frac{2 z t}{1 - \eta z +\sqrt{(1-\eta z)^2-4 z^2\sigma}}
\]
\end{proposition}

\begin{pf} From \rf{Rev Mart} we get $\EE(\EE(\XX_t|\XX_u)^2)=t^2/u$, so applying $\EE(\cdot)$ to \rf{QV'}
we see that $C_{t,u}= t(u-t)/(u+\sigma)$.

Since $\XX_t-\XX_s$ has the same distribution as $\XX_{t-s}$, it suffices to determine the distribution of $\XX_t$.
Fix $t>0$ and let $\XX=\XX_{t}/\sqrt{2t}$ and $\YY=(\XX_{2t}-\XX_t)/\sqrt{2t}$.
Then $\XX,\YY$ are free, centered, and identically distributed. From \rf{Rev Mart}, it
follows that $\EE(\XX|\XX+\YY)=(\XX+\YY)/2$, which implies \rf{LR}.

%Under \rf{LR} we have
%\[\EE((\YY-\XX)^2|\YY+\XX)=\EE((2\XX-(\YY+\XX))^2|\YY+\XX)\]
%\[=4\EE(\XX^2|\XX+\YY)-(\XX+\YY)^2.\]
Assumption \rf{QV'} gives
%\[
%\EE\left.\left(\XX^2\right|\XX+\YY\right)=\frac{1}{4}(\XX+\YY)^2
%+\frac{1}{2+\beta/t}\left(\II+\frac{\alpha}{2\sqrt{t}}(\XX+\YY)+\frac{\beta}{4 t}(\XX+\YY)^2\right).
%\]
\[
\Var(\XX|\XX+\YY)= \frac{1}{4+2\sigma/t}\left(1+\frac{\eta}{\sqrt{2t}}(\XX+\YY)+\frac{\sigma}{2t}(\XX+\YY)^2\right).
\]
Thus \rf{QV} holds with parameters $a=\eta/\sqrt{2t}$,
$b=\sigma/(2t)$. With $\alpha=1/2$, Theorem \ref{T1} says that
 $\XX_t/\sqrt{t}$ has free Meixner $\mu_{\eta/\sqrt{t},\sigma/t}$ law, so
a dilation of \rf{R0} gives the $R$-transform of $\XX_t$.
%
%We observe that the orthogonal polynomials $p_n(x)$ orthogonal with respect to
%the distribution of $\XX_t$
%satisfy the reparametrized version of the recurrence \rf{ccc} which for $n\geq 2$ says
%$$
%xp_n(x)=p_{n+1}(x)+\alpha p_n(x)+(t+\beta) p_{n-1}(x).
%$$

Since the law of $\XX_t$ is a non-negative  measure, from \rf{ccc}
we get $1+\sigma/t\geq 0$. As $t>0$ can be arbitrarily small, we
deduce that $\sigma\geq 0$.
\end{pf}

The next result gives a converse
to \cite[Corollary 7.2]{Capitaine-Casalis04}, and was  inspired by the characterization of Wishart matrices in
\cite[Theorem 4]{Bobecka-Wesolowski02}. Our proof relies on conditional moments and Theorem \ref{T1}.
This method does not work  for Wishart matrices,
see \cite[page 582]{Letac-Massam-98}, who characterize  Wishart matrices by conditional moments of other
 quadratic expressions.

Recall that a random variable $\SSS\in\calA$ is strictly positive if its law is
supported on $[a,b]$ for some $a>0$; in this case $\EE(\SSS)>0$, and from the functional calculus
 the inverse $\SSS^{-1}$ exists and is also strictly positive.
%If $\XX,\YY$ are strictly positive then $\SSS=\XX+\YY$ is strictly positive.
(See \cite[Vol. I]{TakesakiI-III}).

%rownowazne warunki:
%1.x ma spectum dodadnie
%2.jest postaci x=aa* dla pewnego  a  z algebry A.
%3. w pewnej wiernej*- reprezentacji  jego obraz jest operatorem dodatnim.
%
%Co do odwrotnosci to trzeba zalozyc ze element x jest scisle dodatni
%czyli
%spectrum x zawarte w dodatniej prostej i wtedy odwrotny jest tez dodadtni.
%Tutaj mozna zacytowac Takesaki ,3 tomy jego monografii albo Sakai,C* and
%W* algebras.

\begin{proposition}\label{Wishart}
Suppose random variables $\XX,\YY\in\calA$ are non-degenerate,
free,  and such that $\SSS=\XX+\YY$ is strictly positive.
Let $\ZZ=\SSS^{-1/2}\XX\SSS^{-1/2}$.
If $\ZZ$ and $\SSS$ are free, then $\XX$ and $\YY$ have free-Poisson type laws. Moreover, $\EE(\XX)>0$ and the
centered standardized $\XX$ has the free Meixner
 law $\mu_{a,0}$ with $a=\sqrt{\Var(\XX)}/\EE(\XX)$.
\end{proposition}
\begin{pf} We verify that the assumptions of Theorem \ref{T1} are satisfied.
For this proof, we denote the centering operation by $\UU^\circ=\UU-\EE(\UU)\II$.
Denote $\sigma_\XX^2=\Var(\XX), \sigma_\YY^2=\Var(\YY)$, $m_\XX=\EE(\XX),m_\YY=\EE(\YY)$. Since $\SSS$ is strictly positive,
$m_\XX+m_\YY>0$.

We have
$$
\EE(\XX|\SSS)=\EE(\ZZ)\SSS,
$$
as by tracial property and freeness
$\EE(\XX\SSS^n)=\EE(\SSS^{1/2}\ZZ\SSS^{n+1/2})=\EE(\ZZ\SSS^{n+1})=\EE(Z)\EE(\SSS^{n+1})$.
This verifies \rf{LR}.

Applying $\EE()$ we get $\EE(\ZZ)=m_\XX/({m_{\XX}+m_\YY})=\alpha$, so after centering
$$
\EE(\XX^\circ|\SSS^\circ)=\frac{m_\XX}{m_{\XX}+m_\YY}\SSS^\circ.
$$

From \rf{LR} we get
\begin{equation}\label{box}
\sigma_\XX^2=\frac{m_\XX}{m_{\XX}+m_\YY}(\sigma_\XX^2+\sigma_\YY^2).
\end{equation} By non-degeneracy assumption $\sigma_\XX^2>0$,
 this implies that $m_\XX>0$. (By symmetry, $m_\YY>0$, too.)

We now verify that $\Var(\XX|\SSS)$ is a linear function of $\SSS$.
Using tracial property and freeness of $\SSS,\ZZ$, we have
\[\EE(\XX^2\SSS^m)=\EE(\ZZ\SSS\ZZ\SSS^{m+1})=\EE(\ZZ\SSS(\ZZ^\circ+\alpha\II)\SSS^{m+1})
\]
\[
=\alpha\EE(\ZZ\SSS^{m+2})+\EE(\ZZ(\SSS^\circ+(m_{\XX}+m_\YY)\II)\ZZ^\circ\SSS^{m+1})\]
\[
=\alpha^2\EE(\SSS^{m+2})+(m_{\XX}+m_\YY)\EE(\ZZ\ZZ^\circ\SSS^{m+1})+\EE(\ZZ\SSS^\circ\ZZ^\circ\SSS^{m+1})
\]
\[
=\alpha^2\EE(\SSS^{m+2})+(m_{\XX}+m_\YY)\Var(\ZZ)\EE(\SSS^{m+1})+\EE(\ZZ\SSS^\circ\ZZ^\circ\SSS^{m+1}).
\]
Continuing in the same manner, we use freeness to verify that the last term vanishes:
\[\EE(\ZZ\SSS^\circ\ZZ^\circ\SSS^{m+1})=\EE((\ZZ^\circ+\alpha \II)\SSS^\circ\ZZ^\circ\SSS^{m+1})\]\[=
\alpha
\EE(\ZZ^\circ)\EE(\SSS^{m+2})+\EE(\ZZ^\circ\SSS^\circ\ZZ^\circ\SSS^{m+1})\]
\[
=0+\EE(\ZZ^\circ\SSS^\circ\ZZ^\circ)\EE(\SSS^{m+1})=0.
\]
Therefore,
\[\EE(\XX^2\SSS^m)=\EE\left(\left(\alpha^2\SSS^2+(m_{\XX}+m_\YY) \Var(\ZZ)\SSS\right)\SSS^{m}\right),\]
which  by Lemma \ref{Conditionings}(iii) implies that
\[\EE(\XX^2|\SSS)=\alpha^2\SSS^2+(m_{\XX}+m_\YY) \Var(\ZZ)(\SSS^\circ)+(m_{\XX}+m_\YY)^2\Var(\ZZ)\II.\]
Normalizing the variables we get
\[
\Var\left.\left(\frac{1}{\sqrt{\sigma_\XX^2+\sigma_\YY^2}}\XX\right|\SSS\right)=
\frac{(m_{\XX}+m_\YY)^2}{\sigma_\XX^2+\sigma_\YY^2}\Var(\ZZ)
\left(1+\frac{\sqrt{\sigma_\XX^2+\sigma_\YY^2}}{m_{\XX}+m_\YY}\frac{\SSS^\circ}{\sqrt{\sigma_\XX^2+\sigma_\YY^2}}\right).
\]
Thus \rf{QV} holds with $a=\frac{\sqrt{\sigma_\XX^2+\sigma_\YY^2}}{m_{\XX}+m_\YY}$ and $b=0$.

By Theorem \ref{T1}(ii) we see that  $\XX$ is free Poisson type,
and $\XX^\circ/\sigma_\XX$ is free Meixner $\mu_{a/\sqrt{\alpha},0}$
with parameter $a/\sqrt{\alpha}=\frac{\sigma_\XX^2+\sigma_\YY^2}{\sigma_\XX(m_\XX+m_\YY)}=\sigma_\XX/m_\XX$,
see \rf{box}.
This also determines
$
\Var(\ZZ)=\frac{\sigma_\XX^2\sigma_\YY^2}{(\sigma_\XX^2+\sigma_\YY^2)(m_{\XX}+m_\YY)^2}
$.

\end{pf}

Similar reasoning gives the following free analog of Lukacs' theorem \cite{Lukacs-55}.
We do not know whether the property we assume in fact holds for the free gamma law.
\begin{proposition}\label{P free gamma}
Suppose random variables $\XX,\YY\in\calA$ are non-degenerate, i.e.  $\sigma=\sqrt{\Var(\XX)}>0$, free, identically distributed, and strictly positive;
in particular, $m=\EE(\XX)>0$.
For $\SSS=\XX+\YY$ let $\ZZ=\SSS^{-1}\XX^2\SSS^{-1}$.
If $\ZZ$ and $\SSS$ are free, then $\XX$ has free-gamma type law $\mu_{2a,a^2}$ with  $a=\sigma/m$.
\end{proposition}
\begin{pf}
By exchangeability, $\EE(\XX|\SSS)=\SSS/2$, which implies \rf{LR} with $\alpha=1/2$.
By freeness, $\EE(\XX^2|\SSS)=\SSS\EE(\ZZ |\SSS)\SSS=\EE(\ZZ )\SSS^2$. This shows that
\[
\Var(\XX|\SSS)=c\SSS^2,
\]
where $c=\EE(\ZZ )-1/4\geq 0$.
 After centering and normalizing by $\sqrt{2}\sigma$, this implies \rf{QV} with $a/\sqrt{\alpha}=2\sigma/m$, $b/\alpha=\sigma^2/m^2$.
(The latter also determines
$c=\sigma^2/(m^2+2\sigma^2)$.) %Poprawka Wednesday, May 18, 2005 at 00:10
\end{pf}

Next we deduce from  Theorem \ref{T1} a simple variant of \cite[Theorem 5.3]{Nica-96}.
(For related characterizations under more general concepts of non-commutative independence,
see \cite{Hegerfeldt92} and \cite[Proposition 2.5 and Section 4]{Lehner03a}.)

\begin{corollary}[Nica] Suppose random variables $\XX,\YY\in\calA$ are free, non-degenerate, %i.e. $\Var(\XX)\Var(\YY)>0$,
 and such that $\XX+\YY$ and $\XX-\YY$ are free.
Then $\XX$ has the semicircle law.
\end{corollary}
\begin{pf} Changing the random variables to $\XX-\EE(\XX)\II$ and $\YY-\EE(\YY)\II$ preserves the assumptions and the conclusion.
Therefore, without loss of generality we may assume that $\EE(\XX)=\EE(\YY)=0$.
 Since $\XX+\YY$ and $\XX-\YY$ are free, by Lemma \ref{Conditionings}(ii) we have
$\EE(\YY-\XX|{\XX+\YY})=0$, so \rf{LR} holds with $\alpha=1/2$. Moreover,
 $\EE(\XX(\XX+\YY)^n)=\EE(\YY(\XX+\YY)^n)$ for all
$n$, which by \rf{Cumulants} and freeness implies that
\[R_n(\XX)=R_n(\XX,\XX+\YY,\XX+\YY,\dots,\XX+\YY)\]
\[=R_n(\YY,\XX+\YY,\XX+\YY,\dots,\XX+\YY)=R_n(\YY),\]
so $\XX,\YY$ have the same law. Therefore, we can standardize
$\XX,\YY$, dividing both by the same number $\sqrt{\Var(\XX)}>0$.
This operation preserves the freeness of $\XX+\YY$, $\XX-\YY$
and shows that without loss of generality we may assume that
$\XX,\YY$ are standardized with mean $0$ and variance $1$.

Using freeness of $\YY-\XX$ and $\XX+\YY$ again, we get
$\EE((\YY-\XX)^2|{\XX+\YY})=\EE((\YY-\XX)^2)\II=2\II$. Thus
$\Var(\YY|\XX+\YY)=2$ by \rf{Var2} and Theorem \ref{T1} says that
$\XX$ has the Meixner-type law $\mu_{0,0}$, which is the semicircle
law.
\end{pf}

We now  deduce a version of \cite[Corollary 2.4]{Hiwatashi-Masaru-Yoshida99}, with the assumption of
the freeness of the sample mean and the sample variance slightly relaxed.
\begin{corollary}%[Hiwatashi-Masaru-Yoshida]
Let $\XX_1,\XX_2,\dots,\XX_n$ be free identically distributed random
variables, and we put $\overline{\XX}=\frac1n\sum_{j=1}^n\XX_j$
(sample mean), $\VV=\frac1n\sum_{j=1}^n(\XX_j-\overline{\XX})^2$ (sample
variance). If $n\geq 2$ and $ \EE(\VV|\overline{\XX})$ is a multiple
of identity (in particular, if  $\overline{\XX}$ and $\VV$ are
free), then $\XX_1$ has the semicircle law.
\end{corollary}
\begin{pf}
Subtracting $\EE(\XX_1)\II$ from all random variables, without loss
of generality we may assume that $\EE(\XX_1)=0$. If
$\EE(\XX_1^2)=0$, then $\XX_1=0$ has  the (degenerate) semicircle
law. Otherwise, we rescale the random variables, and reduce the
problem to the case $\EE(\XX_1)=0$, $\EE(\XX_1^2)=1$.

We now verify that the assumptions of   Theorem \ref{T1}(i) hold with free random variables
 $\XX=\XX_1/\sqrt{n}$ and $\YY=(\XX_2+\XX_3+\dots+\XX_n)/\sqrt{n}$ .
By exchangeability, compare proof of \rf{Rev Mart}, we have
%$$\EE\left(\XX_1\left|\sum_{j=1}^n\XX_j\right.\right)=\frac1n\sum_{j=1}^n\EE\left(\XX_j\left|\sum_{j=1}^n\XX_j\right.\right)=
%\EE\left(\frac1n\sum_{j=1}^n\XX_j\left|\sum_{j=1}^n\XX_j\right.\right)=\frac1n\sum_{j=1}^n\XX_j.$$
$$\EE\left(\XX_1\left|\overline{\XX}\right.\right)=\frac1n\sum_{j=1}^n\EE\left(\XX_j\left|\sum_{j=1}^n\XX_j\right.\right)=
\EE\left(\frac1n\sum_{j=1}^n\XX_j\left|\sum_{j=1}^n\XX_j\right.\right)=\overline{\XX}.$$

Therefore, \rf{LR} holds with $\alpha=1/n$. Using \rf{Var2}
and  exchangeability we get
 $$\Var(\XX|\XX+\YY)=\EE((\XX_1-\overline{\XX})^2|\overline{\XX})=\frac1n\sum_{j=1}^n\EE((\XX_j-\overline{\XX})^2|\overline{\XX})=
\EE(\VV|\overline{\XX})=C\II,$$
 verifying \rf{QV} with  $a=b=0$.
Since  $\YY$ is non-degenerate for $n\geq 2$, Theorem \ref{T1}(i) implies that
 $\XX_1$  has the semicircle law $\mu_{0,0}=\omega_{0,1}$.
\end{pf}

%We can also recover the CLT:
%\begin{corollary} If $\XX_1,\XX_2,\dots,$ are free identically distributed and $E(X_j)=0$, $E(X_j^2=1$ then
%\end{corollary}

\section{Proof of Theorem \ref{T1}}\label{proof of T1}
%It is easy to see that \rf{LR} implies that $(1-\rho)\YY,\rho\XX$ have the same moments. Indeed,
Since $\SSS:=\XX+\YY$ is non-degenerate, without loss of generality we may assume that
$\EE(\XX)=\EE(\YY)=0$,
$\EE(\XX^2)+\EE(\YY^2)=1$.
Applying $\EE(\cdot)$ to both sides of \rf{LR}  we see that $\alpha_0=0$. Multiplying both sides of \rf{LR} by $\SSS$ and applying $\EE(\cdot)$
we get $\alpha=\EE(\XX^2)$.

Denote  $\beta=\EE(\YY^2)=1-\alpha$ and let
$\VV=\beta\XX-\alpha\YY$. From \rf{LR} it follows that
$\EE(\YY|\SSS)=\beta\SSS$ and \rf{Var2} gives
$\Var(\YY|\SSS)=\EE(\VV^2|\SSS)=\Var(\XX|\SSS)$. Thus the
assumptions are symmetric with respect to $\XX,\YY$, and we only
need  to prove that $\XX$ has a free Meixner type law.

Applying $\EE(\cdot)$ to both sides of \rf{QV}, we get $\EE(\VV^2)=C(1+b)$,
so the normalizing constant is $C=(\beta^2\alpha+\alpha^2\beta)/(1+b)=\alpha\beta/(1+b)$ as claimed.

Next, we establish the following identity:
\begin{equation}\label{P*}
R_n(\XX)=\alpha R_n(\SSS) \mbox{ and } R_n(\YY)=\beta R_n(\SSS).
\end{equation}
We prove this by induction. Since the variables are centered, $R_1(\XX)=R_1(\YY)=R_1(\SSS)=0$. Suppose \rf{P*} holds
true for some $n\geq 1$. Since \rf{LR} implies that $\EE(\XX\SSS^{n})=\alpha\EE(\SSS^{n+1})$, expanding both sides of
this identity into free cumulants we get
\[R_{n+1}(\XX,\SSS,\SSS,\dots,\SSS)+
\sum_{k=1}^n\sum_{\calV=\{B_0,B_1,\dots,B_k\}}R_{B_0}(\XX,\SSS,\SSS,\dots,\SSS)\prod_{j=1}^kR_{|B_j|}(\SSS)
\]
\[
=\alpha R_{n+1}(\SSS)+
\alpha \sum_{k=1}^n\sum_{\calV=\{B_0,B_1,\dots,B_k\}}\prod_{j=0}^kR_{|B_j|}(\SSS).
\]
Then \rf{P*} for $n+1$ follows from induction assumption, as
$R_{n+1}(\XX,\SSS,\SSS,\dots,\SSS)=R_{n+1}(\XX,\XX+\YY,\XX+\YY,\dots,\XX+\YY)=R_{n+1}(\XX)$
by the freeness of $\XX,\YY$.

In particular,
$R_n(\VV,\SSS,\SSS,\dots,\SSS)=R_n(\beta\XX-\alpha\YY,\SSS,\SSS,\dots,\SSS)=\beta R_n(\XX)-\alpha R_n(\YY)$,
 so \rf{P*} implies
\begin{equation}\label{P**}
R_n(\VV,\SSS,\SSS,\dots,\SSS)=0, \; n\geq 1.
\end{equation}
Similarly,  $R_n(\VV,\VV, \SSS,\SSS,\dots,\SSS)=\beta^2R_n(\XX)+\alpha^2R_n(\YY)$, so \rf{P*} implies
\begin{equation}\label{P***}
R_n(\VV,\VV,\SSS,\SSS,\dots,\SSS)=\alpha\beta R_n(\SSS), \; n\geq 2.
\end{equation}

Denote  $m_n=\EE(\SSS^n)$ and let
\[
M(z)=\sum_{n=0}^\infty z^n m_n
\]
be the moment generating function.
\begin{lemma} $M(z)$ satisfies  quadratic equation
\begin{equation} \label{MMM}
(z^2+az+b)M^2-(1+az+2b)M+1+b=0.
\end{equation}
\end{lemma}
\begin{pf}
Multiplying \rf{QV} by $\SSS^n$  for $n\geq 0$ and applying $\EE(\cdot)$ we obtain
\begin{equation}\label{**+}
\EE(\VV^2\SSS^n)=\frac{\alpha\beta}{1+b}( m_n+ a m_{n+1} + b m_{n+2}).
\end{equation}
Expanding the left hand side into the free cumulants we see that
\begin{equation}\label{step1}
\EE(\VV^2\SSS^n)=\sum_{\calV\in\NC(n+2)} R_\calV(\VV,\VV,\SSS,\SSS,\dots,\SSS).
\end{equation}
Since $R_1(\SSS)=\EE(\SSS)=0$, the sum in \rf{step1} can be
restricted to partitions that have no singleton blocks.

Let $\widetilde{\NC}(n+2)$ be the set of all non-crossing partitions of $\{1,2,\dots,n+2\}$
which separate $1$ and $2$ and have no
singleton blocks. Let $\widetilde{\widetilde{\NC}}(n+2)$ denote  the set of
all non-crossing partitions of $\{1,2,\dots,n+2\}$
with the first two elements in the same block  and which have no
singleton blocks.
By \rf{P**}, if a partition $\calV$ separates the first two elements of $\{1,2,\dots,n+2\}$, then
$R_\calV(\VV,\VV,\SSS,\SSS,\dots,\SSS)=0$. Thus the sum in \rf{step1} can be taken over $\widetilde{\widetilde{\NC}}(n+2)$.

If $\calV\in\widetilde{\widetilde{\NC}}(n+2)$ is such a partition,
 then from \rf{P***} we have
\[R_\calV(\VV,\VV,\SSS,\SSS,\dots,\SSS)=\alpha\beta R_\calV(\SSS,\SSS,\dots,\SSS).\]
This shows that  we can eliminate $\VV$  from the right hand side of \rf{step1}. Thus
\[
\EE(\VV^2\SSS^n)=\alpha\beta\sum_{\calV\in\NC(n+2)\setminus \widetilde{\NC}(n+2)}R_\calV(\SSS)\]
\[=
\alpha\beta\sum_{\calV\in\NC(n+2)}R_\calV(\SSS)-\alpha\beta\sum_{\calV\in \widetilde{\NC}(n+2)}R_\calV(\SSS)
=\alpha\beta m_{n+2}-s,
\]
where
\[
s=\alpha\beta\sum_{\calV\in\widetilde{\NC}(n+2)}R_\calV(\SSS).\]

 Since $\widetilde{\NC}(n+2)$ has no singleton blocks, for every
$\calV\in\widetilde{\NC}(n+2)$ there is an index
$k=k(\calV)\in\{3,4,\dots,n+2\}$ such that $k$ is the second
left-most element of the block containing $1$; for example, in the
partition shown in Figure \ref{Fig 1}, $k_\calV=r$. This decomposes
$\widetilde{\NC}(n+2)$ into the $n$ classes
$\widetilde{\NC}_j=\{\calV\in\widetilde{\NC}(n+2): k(\calV)=j+2\}$,
$j=1,2,\dots,n$.

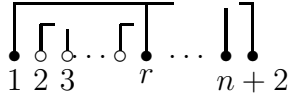
\begin{figure}[tbh]
\begin{picture}(100,110)(-120,-10)

\put(10, 40){\circle*{4}} \put(20, 40){\circle{4}} \put(30,
40){\circle{4}} \put(60, 40){\circle*{4}} \put(50, 40){\circle{4}}
\put(90, 40){\circle*{4}} \put(100, 40){\circle*{4}}
\put(10,40){\makebox(0,-5)[ct]{$1$}}
\put(20,40){\makebox(0,-5)[ct]{$2$}}
\put(30,40){\makebox(0,-5)[ct]{$3$}}
\put(45,40){\makebox(0,0)[rc]{$\dots$}}
\put(75,40){\makebox(0,0)[cc]{$\dots$}}

\put(60,40){\makebox(0,-5)[ct]{$r$}}
\put(100,40){\makebox(0,-5)[ct]{$n+2$}} \thicklines
\put(10,40){\line(0,1){20}} \put(60,40){\line(0,1){20}}
\put(90,40){\line(0,1){18}} \put(100,40){\line(0,1){20}}
\put(100,60){\line(-1,0){5}} \put(20,42){\line(0,1){10}}
\put(30,42){\line(0,1){8}} \put(50,42){\line(0,1){10}}
\put(20,52){\line(1,0){5}}%\put(30,50){\line(1,0){5}}
\put(55,52){\line(-1,0){5}} \put(10,60){\line(1,0){60}}
\end{picture}
\caption{$\calV\in\widetilde{\NC}_{j}$ with $j=r-2$ is decomposed
into two partitions, the first one partitioning the white circles
and the second one partitioning the black circles.\label{Fig 1}}
\end{figure}

Each of the sets $\widetilde{\NC}_j$ is in one-to-one correspondence with the product
\[\NC(j)\times \widetilde{\widetilde{\NC}}(n+2-j).\]
Indeed, the blocks of each partition in $\widetilde{\NC}_j$ consist
of the partition of $\{2,3,\dots,j+1\}$ which can be uniquely
identified with the appropriate partition in $\NC(j)$, and the
remaining blocks which  partition the $(n+2-j)$-element set
 $\{1,j+2,j+3,\dots,n+2\}$ under the additional constraint that the first two elements $1,j+2$
  are in the same block, see Figure \ref{Fig 1}.
 These remaining blocks can therefore be uniquely identified with the partition in
 $\widetilde{\widetilde{\NC}}(n+2-j)$. This gives
 \[
 s=\alpha\beta\sum_{j=1}^n \sum_{\calV\in \widetilde{\NC}_j}R_\calV(\SSS)=\sum_{j=1}^n\sum_{\calV\in\NC(j)}R_\calV(\SSS)
 \sum_{\calV\in\widetilde{\widetilde{\NC}}(n+2-j)}\alpha\beta R_\calV(\SSS)\]\[
=\sum_{j=1}^n m_j \EE(\VV^2\SSS^{n-j}),
\]
which establishes the identity
\[\EE(\VV^2\SSS^n)=\alpha\beta m_{n+2}-\sum_{j=1}^n m_j \EE(\VV^2\SSS^{n-j}).
\]
Since $m_0=1$ and $\alpha\beta>0$ as $\XX,\YY$ are non-degenerate, combining the above formula with \rf{**+} we get
\begin{equation}\label{***}
m_{n+2}=\frac{1}{1+b}\sum_{j=0}^n m_j(m_{n-j}+am_{n+1-j}+bm_{n+2-j}).
\end{equation}
Using the fact that $m_1=0$, from \rf{***} we obtain
\[
M(z)-1=\frac{z^2}{1+b}M^2(z)+\frac{az}{1+b}M(z)(M(z)-1)+\frac{b}{1+b}M(z)(M(z)-1),
\]
which is equivalent to \rf{MMM}.
\end{pf}
From \rf{MMM} we see that
\[
M(z)=\frac{1+2b+az-\sqrt{(1-az)^2-4z^2(1+b)}}{2(z^2+az+b)}.
\]
From this we calculate the corresponding Cauchy-Stieltjes transform $G_{\SSS}(z)=M(1/z)/z$
and the $R$-transform of $\SSS$ which for $b\ne0$ takes the form
\[r_{\SSS}(z)=\frac{1-az-\sqrt{(1-az)^2-4z^2b}}{2 zb}.\]
Thus the distribution of $\SSS$ is the free Meixner measure $\mu_{a,b}$ given by \rf{R0}.
The $R$-transform of $\XX$ is
\[r_{\XX}(z)=\alpha\frac{1-az-\sqrt{(1-az)^2-4z^2b}}{2 zb}.\]
see \rf{P*}. After standardization,
\[r_{\XX/\sqrt{\alpha}}(z)=\frac{1-az/\sqrt{\alpha}-\sqrt{(1-az/\sqrt{\alpha})^2-4z^2b/\alpha}}{2 zb/\alpha},\]
so $\XX$ has the free Meixner type law
$\mu_{a/\sqrt{\alpha},b/\alpha}$.

To end the proof, we notice that \rf{ccc} applied to the law of
$\XX$ implies $b/\alpha\geq -1$, see \cite[page 21]{Chihara}. Since
the distribution of $\YY$ must also be of free Meixner type and well
defined, we get $b/\beta\geq -1$. Thus $b\geq
-\min\{\alpha,1-\alpha\}$ as claimed.

\section{Remarks}
\begin{remark}
If $\XX$ is free Poisson with mean $m>0$ then
$\sigma^2=m$, thus $a=1/\sqrt{m}$ and the law of $\XX$ is  $\delta_m\boxplus D_{\sqrt{m}}(\mu_{1/\sqrt{m},0})$.  By
 \cite[Corollary 7.2]{Capitaine-Casalis04}, when $m>1$ the random variables $\ZZ,\SSS$ are indeed free as
 assumed in Proposition \ref{Wishart}.
 But Proposition \ref{Wishart} is not a characterization of all free-Poisson type laws as
 the free-Poisson laws with $m\leq 1$ fail to be strictly positive, and when $m<1$ have an atom at $0$.
\end{remark}
\begin{remark}
Proposition \ref{Wishart} extends with the same proof to the case when
$\SSS=\XX+\YY$ and $\ZZ=\SSS_1^{-1}\XX\SSS_2^{-1}$ are free for any  decomposition of $\SSS=\SSS_1\SSS_2$,
which is the setting of the original Olkin-Rubin \cite{Olkin-Rubin-62} characterization of the Wishart matrices.

The fact that \rf{QV} should hold true with $b=0$ is to be expected from the expression for
 the conditional moment of the square of a Wishart matrix given
 after \cite[Corollary 2.3]{Letac-Massam-98}. \end{remark}

%\begin{remark} \cite[Theorem 2.3]{Hiwatashi-Masaru-Yoshida99} suggests that
%a version of Theorem \ref{T1} might also hold for more general linear and quadratic forms
%in freely identically distributed random variables.
%\end{remark}
\begin{remark}  Our proof of Theorem \ref{T1} does not rely to a significant degree on $*$-operation
 of $\calA$ and could apply to any abstract probability space, see  \cite{Haagerup97}.
One exception is the argument that shows $b\geq -\min\{\alpha,1-\alpha\}$; the setting of von Neumann algebras helps with
 conditional expectations, existence of which would have to be assumed in the more general setting.

In the tracial von Neumann setting, the proof seems to cover free random variables  $\XX,\YY\in\bigcap_{p>1}L_p(\calA,\tau)$.
It would be nice to extend Theorem \ref{T1} to a complete analog of the classical setting, with
random variables in $L_2(\calA,\tau)$ only.
\end{remark}

\begin{remark}[$q$-interpolation]\label{q-interpolation}
For $-1<q\leq 1$ consider the following recurrence
\begin{equation}\label{CCC}
R_{n+1}=a R_n+ b \sum_{j=2}^{n-1}\left[ \begin{array}{c}n-1 \\ j-1 \end{array}\right]_q R_j R_{n+1-j}, \, n\geq 2,
\end{equation}
with the initial values $R_1=0, R_2=1$.
Here we use the standard notation
%\begin{eqnarray*}
%{[n]_{q}} &=&1+q+\dots +q^{n-1},
%{[n]_{q}!} &=&[1]_{q}[2]_{q}\dots [n]_{q},
%\\ \left[ \begin{array}{c}n \\ k \end{array}\right] _{q} &=&\frac{[n]_{q}!}{[n-k]_{q}![k]_{q}!},
%\end{eqnarray*}%
\[
{[n]_{q}} =1+q+\dots +q^{n-1}, \;
{[n]_{q}!} =[1]_{q}[2]_{q}\dots [n]_{q}, \;
\\ \left[ \begin{array}{c}n \\ k \end{array}\right] _{q} =\frac{[n]_{q}!}{[n-k]_{q}![k]_{q}!},
\]
with the usual conventions $[0]_{q}=0,[0]_{q}!=1$.

The proof of Theorem \ref{Laha-Lukacs} relies on the differential equation for the derivative of the
logarithm of the characteristic function, i.e., for the $R_1$-transform of Nica  \cite{Nica-95}.
This differential equation is equivalent to \rf{CCC}  holding with $q=1$ for the classical cumulants $c_n(X+Y)$.

 One can check that \rf{R0} is equivalent to the quadratic equation
\[
zb r^2-(1-az)r+z=0
\] for the $R$ transform $r=r(z)$.
Since $r(z)=\sum_{k=0}^\infty R_{k+1}(\mu)z^k$, this implies that when  $q=0$,
the recurrence \rf{CCC} holds  for the free cumulants $R_n(\XX+\YY)$.
For a related observation see \cite[Propositions 1 and 2]{Anshelevich-04}.
  \end{remark}

%\begin{remark}
%%Wesolowski \cite[Theorem 3.2 and Eqn (*)]{Wesolowski89a}  derives another recurrence for the
%%moments in a related setting.
%Formula \rf{***} is a $q=0$ version of the identity
% $$
%m_{n+2}=
% .$$
% \fbox{Jacka dowod post factum daje jakas rekurencje!}
%  The $q=1$ version for the classical situation is  implicit in
%\cite[Theorem 3.2]{Wesolowski89a}.
%\end{remark}
\begin{remark}Recall that a classical version of a non-commutative process $(\XX_t)$ is a classical process $(X_t)$ on some probability space $(\Omega,\calF,P)$
such that
\[
\EE(\XX_{t_1}\XX_{t_2}\dots\XX_{t_n})=E(X_{t_1}X_{t_2}\dots X_{t_n})
\]
for all $0\leq t_1\leq t_2\leq \dots \leq t_n$. In \cite[Section 4]{BKS97},
 the authors show that classical versions exist for all  $q$-Gaussian Markov processes.
From \cite[Theorem 4.2]{Biane98} it follows that the classical
version of the process $(\XX_t)$ from  our Proposition \ref{SP}
exists and is a classical Markov process. This is the same
process  that appears in \cite[Theorem 3.5]{Bryc-Wesolowski-03}
when the parameters are $q=0$, $\theta=\eta$, $\tau=\sigma$.
It is interesting to note that the conditional variances of the classical versions are also quadratic, % but different from \rf{QV},
see
\cite[Theorem 4.3]{Bryc-Wesolowski-03} and that in
the classical case there is a family of Markov processes with the laws that for $0\leq q\leq 1$
interpolate between the free Meixner laws of Theorem \ref{T1} and the classical Meixner laws of Theorem \ref{Laha-Lukacs}.
\end{remark}

\begin{remark}
The free Meixner laws with $-1\leq b<0$  are not $\boxplus$-infinitely divisible (\cite{Saitoh-Yoshida01}), but they are
infinitely divisible with respect to the $c$-convolution
\cite{Bozejko-Leinert-Speicher}, and appear in generalized limit theorems
\cite{Bozejko-Wysoczanski01}.
 \end{remark}
\begin{remark} The Catalan numbers show up as cumulants of the free Meixner laws in several different situations.
Firstly, for the symmetric free Meixner law $\mu_{0,b}$, the cumulants are
$R_{2k+1}=0$ and $R_{2k+2}=\frac{1}{k+1}\left(^{2k}_{k}\right)b^{k}$, $k\geq 0$; in particular, the Catalan numbers appear
as cumulants of the $\boxplus$-infinitely divisible law $\mu_{0,1}$ and, with alternating signs, as cumulants of
 the two-point law $\mu_{0,-1}$.
Secondly,  the cumulants of the standardized free Gamma type law $\mu_{2a,a^2}$ are $R_1=0$ and
$R_{k+1}=\frac{1}{k+1}\left(^{2k}_k\right)a^{k-1}$, $k\geq 1$.
Compare the Delaney triangle in \cite{Bozejko-Wysoczanski01}.
\end{remark}
%\begin{pf} The first statement is trivial. The second, is a computer observation.
%\end{pf}
\begin{remark}
It would be interesting to know whether Theorem \ref{T1} admits random matrix models, i.e. whether
there are pairs of independent random  matrices $\mX,\mY$ which have the same
law that is invariant under orthogonal transformations, i.e.
 $U\mX U^T$ has the same law as $\mX$ for any deterministic orthogonal matrix $U$, satisfy
$E(\mX^2|\mX+\mY)= C (4 + 2a (\mX+\mY) + b^2(\mX+\mY)^2)$ and are asymptotically free.

All $\boxplus$-infinitely divisible
laws have matrix models, see
\cite{Benaych-Georges04} and \cite{Cabanal-Duvillard04}, see also \cite[Section 4.4]{Hiai-Petz},
 but it is not clear whether one can preserve the
quadratic form of the conditional variance.
\end{remark}
\begin{remark}
Regarding Proposition \ref{P free gamma}, it would be interesting to know whether there are i.i.d.
symmetric random matrices  $\mX,\mY$ with independent
$\mS=\mX+\mY$ and $\mZ=\mS^{-1}\mX^2\mS^{-1}$ and with the law of $\mX$ that is invariant under orthogonal transformations.
Wishart $n\times n$ matrices with scale parameter $\mI$
cannot have the above property for large $n$, as their asymptotic distribution is different.

{\bf Note added late:}  According to G. L\'etac \cite{Letac05},
 positive-definite $n\times n$ matrices with the above property do not exist for $n>1$.
\end{remark}
\subsection*{Acknowledgement}  We thank M. Anshelevich, G. Letac, and J. Weso{\l}owski for helpful comments and references.
We thank the anonymous referee for very careful reading of the submitted manuscript.
The first author would like to thank
for fantastic working conditions at the Department of Mathematical Sciences of the University of Cincinnati during
his visit to Cincinnati in September-October 2004. The second author would like to thank
M. Anshelevich for raising the question of $q$-generalizations of \cite{Laha-Lukacs60}.

%\tolerance=70000
%%BibTeX
%%\bibliographystyle{apalike}
%\bibliographystyle{plain}
%%\bibliographystyle{elsart-num}
%\bibliography{Vita,Wesol,free-lalu,q-reg}

\begin{thebibliography}{10}

\bibitem{Accardi-Bozejko98}
Luigi Accardi and Marek Bo{\.z}ejko.
\newblock Interacting {F}ock spaces and {G}aussianization of probability
  measures.
\newblock {\em Infin. Dimens. Anal. Quantum Probab. Relat. Top.},
  1(4):663--670, 1998.

\bibitem{Anshelevich01a}
Michael Anshelevich.
\newblock Partition-dependent stochastic measures and {$q$}-deformed cumulants.
\newblock {\em Doc. Math.}, 6:343--384 (electronic), 2001.

\bibitem{Anshelevich01}
Michael Anshelevich.
\newblock Free martingale polynomials.
\newblock {\em Journal of Functional Analysis}, 201:228--261, 2003.
\newblock arXiv:math.CO/0112194.

\bibitem{Anshelevich-04}
Michael Anshelevich.
\newblock {Orthogonal polynomials with a resolvent-type generating function},
  2004.
\newblock arXiv:math.CO/0410482 (to appear in Trans Amer Math Soc.

\bibitem{Barndorff-Nielsen-Thorbj02c}
Ole~E. Barndorff-Nielsen and Steen Thorbj{\o}rnsen.
\newblock L\'evy processes in free probability.
\newblock {\em Proc. Natl. Acad. Sci. USA}, 99(26):16576--16580 (electronic),
  2002.

\bibitem{Barndorff-Nielsen-Thorbj02a}
Ole~E. Barndorff-Nielsen and Steen Thorbj{\o}rnsen.
\newblock Self-decomposability and {L}\'evy processes in free probability.
\newblock {\em Bernoulli}, 8(3):323--366, 2002.

\bibitem{Benaych-Georges04}
Florent Benaych-Georges.
\newblock Classical and free infinitely divisible distributions and random
  matrices.
\newblock Preprint: arXiv:math.PR/0406082.

\bibitem{Berkovici-Pata99}
Hari Bercovici and Vittorino Pata.
\newblock Stable laws and domains of attraction in free probability theory.
\newblock {\em Ann. of Math. (2)}, 149(3):1023--1060, 1999.
\newblock With an appendix by Philippe Biane.

\bibitem{Berkovici-Voiculescu92}
Hari Bercovici and Dan Voiculescu.
\newblock L\'evy-{H}in\v cin type theorems for multiplicative and additive free
  convolution.
\newblock {\em Pacific J. Math.}, 153(2):217--248, 1992.

\bibitem{Biane98}
Philippe Biane.
\newblock Processes with free increments.
\newblock {\em Math. Z.}, 227(1):143--174, 1998.

\bibitem{Bobecka-Wesolowski02}
Konstancja Bobecka and Jacek Weso{\l}owski.
\newblock The {L}ukacs-{O}lkin-{R}ubin theorem without invariance of the
  ``quotient''.
\newblock {\em Studia Math.}, 152(2):147--160, 2002.

\bibitem{BKS97}
Marek Bo{\.z}ejko, Burkhard K{\"u}mmerer, and Roland Speicher.
\newblock $q$-{G}aussian processes: non-commutative and classical aspects.
\newblock {\em Comm. Math. Phys.}, 185(1):129--154, 1997.

\bibitem{Bozejko-Leinert-Speicher}
Marek Bo{\.z}ejko, Michael Leinert, and Roland Speicher.
\newblock Convolution and limit theorems for conditionally free random
  variables.
\newblock {\em Pacific J. Math.}, 175(2):357--388, 1996.

\bibitem{Bozejko-Speicher91}
Marek Bo{\.z}ejko and Roland Speicher.
\newblock {$\psi$}-independent and symmetrized white noises.
\newblock In {\em Quantum probability \& related topics}, QP-PQ, VI, pages
  219--236. World Sci. Publishing, River Edge, NJ, 1991.

\bibitem{Bozejko-Wysoczanski01}
Marek Bo{\.z}ejko and Janusz Wysocza{\'n}ski.
\newblock Remarks on {$t$}--transformations of measures and convolutions.
\newblock {\em Ann. Inst. H. Poincar\'e Probab. Statist.}, 37(6):737--761,
  2001.

\bibitem{Bryc-Wesolowski-03}
W{\l}odzimierz Bryc and Jacek Weso{\l}owski.
\newblock Conditional moments of $q$--{M}eixner processes.
\newblock {\em Probability Theory Related Fields}, 131:415--441, 2005.
\newblock arxiv.org/abs/math.PR/0403016.

\bibitem{Cabanal-Duvillard04}
Thierry Cabanal-Duvillard.
\newblock A matrix repersentation of the {B}ercovici--{P}ata bijection.
\newblock {\em Electronic Journal of Probability}, 10:632--661, 2005.
\newblock Paper 18.

\bibitem{Capitaine-Casalis04}
M.~Capitaine and M.~Casalis.
\newblock Asymptotic freeness by generalized moments for {G}aussian and
  {W}ishart matrices. {A}pplication to beta random matrices.
\newblock {\em Indiana Univ. Math. J.}, 53(2):397--431, 2004.

\bibitem{Chihara}
T.~S. Chihara.
\newblock {\em An introduction to orthogonal polynomials}.
\newblock Gordon and Breach, New York, 1978.

\bibitem{Collins04}
Benoit Collins.
\newblock {Product of random projections, {J}acobi ensembles and universality
  problems arising from free probability}.
\newblock arXiv:math.PR/0406560.

\bibitem{Haagerup97}
Uffe Haagerup.
\newblock On {V}oiculescu's {$R$}-- and {$S$}--transforms for free non-commuting
  random variables.
\newblock In {\em Free probability theory (Waterloo, ON, 1995)}, volume~12 of
  {\em Fields Inst. Commun.}, pages 127--148. Amer. Math. Soc., Providence, RI,
  1997.

\bibitem{Haagerup-Larsen}
Uffe Haagerup and Flemming Larsen.
\newblock Brown's spectral distribution measure for {$R$}--diagonal elements in
  finite von {N}eumann algebras.
\newblock {\em J. Funct. Anal.}, 176(2):331--367, 2000.

\bibitem{Hegerfeldt92}
Gerhard~C. Hegerfeldt.
\newblock A quantum characterization of {G}aussianness.
\newblock In {\em Quantum probability \& related topics}, QP--PQ, VII, pages
  165--173. World Sci. Publishing, River Edge, NJ, 1992.

\bibitem{Hiai-Petz}
Fumio Hiai and D{\'e}nes Petz.
\newblock {\em The semicircle law, free random variables and entropy},
  volume~77 of {\em Mathematical Surveys and Monographs}.
\newblock American Mathematical Society, Providence, RI, 2000.

\bibitem{Hiwatashi-Masaru-Yoshida99}
Osamu Hiwatashi, Masaru Nagisa, and Hiroaki Yoshida.
\newblock The characterizations of a semicircle law by the certain freeness in
  a {$C\sp *$}--probability space.
\newblock {\em Probab. Theory Related Fields}, 113(1):115--133, 1999.

\bibitem{Kesten}
Harry Kesten.
\newblock Symmetric random walks on groups.
\newblock {\em Trans. Amer. Math. Soc.}, 92:336--354, 1959.

\bibitem{Laha-Lukacs60}
R.~G. Laha and E.~Lukacs.
\newblock On a problem connected with quadratic regression.
\newblock {\em Biometrika}, 47(300):335--343, 1960.

\bibitem{Lehner03a}
Franz Lehner.
\newblock Cumulants in noncommutative probability theory. {II}. {G}eneralized
  {G}aussian random variables.
\newblock {\em Probab. Theory Related Fields}, 127(3):407--422, 2003.

\bibitem{Letac05}
G{\'e}rard Letac.
\newblock Remarks on the {B}ozejko and {B}ryc problem.
\newblock Private communication, \LaTeX\ file: July 13, 2005.

\bibitem{Letac-Massam-98}
G{\'e}rard Letac and H{\'e}l{\`e}ne Massam.
\newblock Quadratic and inverse regressions for {W}ishart distributions.
\newblock {\em Ann. Statist.}, 26(2):573--595, 1998.

\bibitem{Lukacs-55}
Eugene Lukacs.
\newblock A characterization of the gamma distribution.
\newblock {\em Ann. Math. Statist.}, 26:319--324, 1955.

\bibitem{Marchenko-Pastur67}
V.A. Marchenko and L.A. Pastur.
\newblock Distribution of eigenvalues for some sets of random matrices.
\newblock {\em USSR, Sb.}, 1:457–--483, 1967.

\bibitem{McKay}
Brendan~D. McKay.
\newblock The expected eigenvalue distribution of a large regular graph.
\newblock {\em Linear Algebra Appl.}, 40:203--216, 1981.

\bibitem{Meixner-40}
J.~Meixner.
\newblock {O}rthogonale {P}olynomsysteme mit einer besonderen {G}estalt der
  erzeugenden {F}unktion.
\newblock {\em Journal of the London Mathematical Society}, 9:6--13, 1934.

\bibitem{Morris82}
Carl~N. Morris.
\newblock Natural exponential families with quadratic variance functions.
\newblock {\em Ann. Statist.}, 10(1):65--80, 1982.

\bibitem{Nica-95}
Alexandru Nica.
\newblock A one-parameter family of transforms, linearizing convolution laws
  for probability distributions.
\newblock {\em Comm. Math. Phys.}, 168(1):187--207, 1995.

\bibitem{Nica-96}
Alexandru Nica.
\newblock {$R$}--transforms of free joint distributions and non-crossing
  partitions.
\newblock {\em J. Funct. Anal.}, 135(2):271--296, 1996.

\bibitem{Nica-Speicher}
Alexandru Nica and Roland Speicher.
\newblock On the multiplication of free {$N$}-tuples of noncommutative random
  variables.
\newblock {\em Amer. J. Math.}, 118(4):799--837, 1996.

\bibitem{Olkin-Rubin-62}
Ingram Olkin and Herman Rubin.
\newblock A characterization of the {W}ishart distribution.
\newblock {\em Ann. Math. Statist.}, 33:1272--1280, 1962.

\bibitem{Saitoh-Yoshida01}
Naoko Saitoh and Hiroaki Yoshida.
\newblock The infinite divisibility and orthogonal polynomials with a constant
  recursion formula in free probability theory.
\newblock {\em Probab. Math. Statist.}, 21(1):159--170, 2001.

\bibitem{Schoutens00}
Wim Schoutens.
\newblock {\em Stochastic processes and orthogonal polynomials}, volume 146 of
  {\em Lecture Notes in Statistics}.
\newblock Springer-Verlag, New York, 2000.

\bibitem{Speicher-97}
Roland Speicher.
\newblock Free probability theory and non-crossing partitions.
\newblock {\em S\'em. Lothar. Combin.}, 39:Art.\ B39c, 38 pp.\ (electronic),
  1997.

\bibitem{TakesakiI-III}
M.~Takesaki.
\newblock {\em Theory of operator algebras. {I-III}}, volume 124 of {\em
  Encyclopaedia of Mathematical Sciences}.
\newblock Springer-Verlag, Berlin, 2002.
\newblock Reprint of the first (1979) edition, Operator Algebras and
  Non-commutative Geometry, 5, 6, 8.

\bibitem{Takesaki}
Masamichi Takesaki.
\newblock Conditional expectations in von {N}eumann algebras.
\newblock {\em J. Functional Analysis}, 9:306--321, 1972.

\bibitem{VDN92}
D.~V. Voiculescu, K.~J. Dykema, and A.~Nica.
\newblock {\em Free random variables}.
\newblock American Mathematical Society, Providence, RI, 1992.

\bibitem{Voiculescu00}
Dan Voiculescu.
\newblock Lectures on free probability theory.
\newblock In {\em Lectures on probability theory and statistics (Saint-Flour,
  1998)}, volume 1738 of {\em Lecture Notes in Math.}, pages 279--349.
  Springer, Berlin, 2000.

\bibitem{Wesolowski89a}
Jacek Weso{\l}owski.
\newblock Characterizations of some processes by properties of conditional
  moments.
\newblock {\em Demonstratio Math.}, 22(2):537--556, 1989.

\bibitem{Wesolowski93}
Jacek Weso{\l}owski.
\newblock Stochastic processes with linear conditional expectation and
  quadratic conditional variance.
\newblock {\em Probab. Math. Statist.}, 14:33--44, 1993.

\bibitem{Wigner58}
Eugene~P. Wigner.
\newblock On the distribution of the roots of certain symmetric matrices.
\newblock {\em Ann. of Math. (2)}, 67:325--327, 1958.

\end{thebibliography}

\end{document}